\documentclass{amsart}

\usepackage[english]{babel}
\usepackage{leftidx}
\usepackage[numbers]{natbib}
\usepackage{bigints}
\usepackage{empheq}
\usepackage{array}
\usepackage{graphicx}
\usepackage{lipsum}
\usepackage{amssymb}
\usepackage{amsthm}
\usepackage{needspace}
\usepackage[makeroom]{cancel}
\newenvironment{myproof}[2] {\paragraph{\textbf{Proof of {#1} {#2} :}}}{\hfill$\square$}
\newtheorem{theorem}{Theorem}[section]
\newtheorem{proposition}[theorem]{Proposition}
\newtheorem{lema}[theorem]{Lemma}
\newtheorem{question-non}[]{}

\newtheorem{cor}[theorem]{Corollary}
\newtheorem{definition}[theorem]{Definition}
\newtheorem{example}[theorem]{Example}

\title{gradient estimates on warped product gradient almost ricci solitons}

\author{Tokura, W. $^{1}$}
\address{$^{1}$ Universidade Federal de Goi\'as, IME, 131, 74001-970, Goi\^ania, GO, Brazil.}
\email{williamisaotokura@hotmail.com $^{1}$}

\author{Adriano, L. $^{2}$}
\address{$^{2}$ Universidade Federal de Goi\'as, IME, 131, 74001-970, Goi\^ania, GO, Brazil.}
\email{levi@ufg.br $^{2}$}

\author{Pina, R. $^{3}$}
\address{$^{3}$ Universidade Federal de Goi\'as, IME, 131, 74001-970, Goi\^ania, GO, Brazil.}
\email{romildo@ufg.br $^{3}$}

\author{Barboza, M. $^{4}$}
\address{$^{4}$ Insituto Federal Goiano, 75790-000, Rodovia Geraldo Silva Nascimento Km 2,5, Uruta\'i, GO, Brazil.}
\email{marcelo.barboza@ifgoiano.edu.br $^{4}$}

\thanks{$^{1,4}$ Supported by CAPES}

\keywords{Warped product, Gradient almost Ricci solitons, Gradient estimates, Einstein manifolds, Quasi-Einstein manifold.}

\subjclass[2010]{53C21, 53C50, 53C25} 

\begin{document}

\begin{abstract}
	In this paper, by slightly modifying Li-Yau's technique so that we can handle drifting Laplacians, we were able to find three different gradient estimates for the warping function, one for each sign of the Einstein constant of the fiber manifold. As an application, we exhibit a nonexistence theorem for gradient almost Ricci solitons possessing certain metric properties on the base of the warped product.
\end{abstract}

\maketitle

\section{Introduction and main results}
\label{intro}

A gradient almost Ricci soliton \cite{pigola2010ricci} is a complete Riemannian manifold $(M^{n}, g)$ together with functions $ h,\rho:M\to\mathbb{R} $ satisfying the equation
\begin{equation}\label{1}Ric_{g}+Hess_{g}h=\rho g,
\end{equation}
where $Ric_{g}$ and $Hess_{g}h$ stand, respectively, for the Ricci tensor and the Hessian of $h$. The function $h:M\to\mathbb{R}$ is usually referred to as the potential function of the soliton in the literature. The quadruple $(M^{n},g,h,\rho)$ is classified into three types according to the sign of $\rho$: expanding if $\rho<0$, steady if $\rho=0$ and shrinking if $\rho>0$. If $\rho$ occurs as a constant, the soliton is usually referred to as a \textit{gradient Ricci soliton}.

From the above it is readily seen that Einstein manifolds, that is, those for which
\[Ric_{g}=\rho g,\]
trivially fulfill the requirement for a gradient almost Ricci soliton with a constant $\rho$ and Killing vector field $\nabla h$, not to exclude the possibility of an also constant $ h $.

Recent works like \cite{de2017gradient, feitosa2015gradient, feitosa2017construction, kim2013warped, pina2017invariant} exhibit the fact that warped product manifolds constitute a highly profitable ground for the construction of gradient almost Ricci solitons. It's worth noticing that [12] contain infinitely many solutions to (1) in closed form.

\begin{definition} Let $(B^{n}, g_{B})$ and $(F^{k}, g_{F})$ be two Riemannian manifolds and $f>0$ on $B$. The product manifold $B\times F$ furnished with metric tensor
	\begin{equation}\label{warped metric}
	g=\pi^{\ast}g_{B}+(f\circ\pi)^{2}\sigma^{\ast}g_{F},
	\end{equation}
	where $\pi:B\times F\rightarrow B$, $\sigma:B\times F\rightarrow F$ are the projections on the first and second factor, respectively, is called \textit{warped product manifold} and will be denote by $B\times_{f}F$. The function $f$ is called \textit{warping function}, $B$ is called the \textit{base} and $F$ the \textit{fiber}.
\end{definition}

In \cite{de2017gradient} the authors investigated warped product gradient Ricci solitons and proved that either the warping function $f$ is constant or the potential function satisfies
\begin{equation*}
 h=h_{B}\circ\pi,\hspace{0,5cm} h_{B}\in C^{\infty}(B).
\end{equation*}

In this very spirit the authors of \cite{feitosa2015gradient} studied warped product gradient almost Ricci solitons with the assumption that
\begin{equation}\label{3'}
\rho=\rho_{B}\circ\pi, \hspace{0,5cm} h=h_{B}\circ\pi,\hspace{0,7cm}\rho_{B},\hspace{.1cm} h_{B}\in C^{\infty}(B),
\end{equation}
and prove the following result:

\begin{proposition}\label{lemma1} Let
	$(B^{n}\times_{f}F^{k},g,h,\rho)$ be a \textbf{gradient almost Ricci soliton}, then the functions $f$, $h_{B}$ defined on $(B^{n},g_{B})$ satisfy
	\begin{equation}\label{111}
	Ric_{g_B}+Hess_{g_B}h_{B}=\rho_{B} g_{B}+\frac{k}{f}Hess_{g_B}f,
	\end{equation}
	and the fiber $(F^{k}, g_{F})$ is \textbf{Einstein} with $Ric_{g_{F}}=\theta g_{F}$, where
	
	\begin{equation}\label{nazareno1}
	\theta=\rho_{B} f^{2}+f\Delta f+(k-1)|\nabla f|^{2}-f\nabla h_{B}(f).
	\end{equation}
\end{proposition}

The authors make use of the equation \eqref{nazareno1} in combination with the maximum principle to provide triviality results for warping functions reaching a maximum. As a consequence one sees that it is impossible to achieve the structure of a gradient almost Ricci soliton on warped product manifolds possessing a compact base and $\rho_{B}f^2\geq\theta$.

Proposition \ref{lemma1} is known to be true for Einstein manifolds $(B^{n}\times_{f}F^{k},g)$ since \cite{kim2013warped}. Besides, with the change $f =\exp{(k^{-1}u)}$ in \eqref{111}, we get that
\begin{equation}\label{222}
Ric_{g_{B}}+Hess_{g_{B}}u-\frac{1}{k}du\otimes du=\rho_{B} g_{B}.
\end{equation}
Objects satisfying equation \eqref{222} are called \textit{quasi-Einstein manifold} \cite{case2011rigidity} and play an important role in the analysis of solutions for \eqref{nazareno1}. For instance, in \cite{case2010nonexistence}, by starting off with \eqref{222}, the author developed a gradient estimate for $u$ and concluded that there are no nonconstant solutions for \eqref{nazareno1} if $\rho_{B}=0$ and $\theta\leq0$. Therefore, there are no nontrivial Ricci flat warped products such that its fiber manifold has a nonpositive Einstein constant. In the same way, in \cite{mastrolia2014some} it's shown some triviality results for the warping function $f$ in terms of gradient estimates, however, in the case $\rho_{B}\leq0$.

By the above exposed it is interesting to investigate whether or not we have local gradient estimates for warping solutions of \eqref{nazareno1} in the general case. By the change $f=v^{\frac{1}{k}}$ in equation \eqref{nazareno1} we obtain that
\begin{equation}\label{9'}
\Delta_{h_{B}}v+\rho_{B} kv-\theta k v^{1-\frac{1}{k}}=0.
\end{equation}
where $\Delta_{w}=e^{w}div(e^{-w}\nabla w)$ is the so called drifting Laplacian on the Bakry-\'Emery geometry. In order to accomplish such a task, we focus our attention on gradient estimates for the positive solutions of the nonlinear equation \eqref{9'} that, in compliment to \cite{feitosa2015gradient}, are defined on a noncompact base. We mainly follow the lines of P. Li and S.T. Yau's proof in \cite{li1986parabolic}.

\begin{theorem}\label{gradient}Let $(B^{n}\times_{f}F^{k},g, h,\rho)$ be a complete gradient almost Ricci soliton with noncompact base satisfying
	$$Ric_{g_B}+Hess_{g_B}h_{B}-\frac{1}{m}dh_{B}\otimes dh_{B}\geq -K,\qquad\Delta_{h_{B}}\rho_{B}\geq0,\qquad|\nabla\rho_{B}|\leq\gamma(2R),$$
	for $K\geq0$ in the metric ball $B(p,2R)$ of the base. Then, for any $\beta\in(0,1)$ with $\beta>1-\frac{2}{k}$, the warping function $f$ satisfies the following gradient estimates:
	
	\begin{itemize}
		\item If $\theta<0$, we have
		
		$$\beta\frac{|\nabla f|^{2}}{f^{2}}+\frac{\rho_B}{k}-\frac{\theta}{kf^2}\leq \frac{(n+m)}{ k^{2}\beta}P+\sqrt{\frac{n+m}{2\beta k^{4}}}Q^{\frac{1}{2}}, \quad \text{in}\quad B(p,R),\hspace{2,5cm}$$
		
		\item If $\theta\geq0$, assume that $f$ is bounded in $B(p,2R)$, then we have
		
		$$\beta\frac{|\nabla f|^{2}}{f^{2}}+\frac{\rho_B}{k}-\frac{\theta}{kf^2}\leq \frac{(n+m)}{k^{2}\beta}\left[P+2\theta M\right]+\sqrt{\frac{n+m}{2\beta k^{4}}}(Q+S)^{\frac{1}{2}}, \quad \text{in}\quad B(p,R),$$
		
	\end{itemize}
	where
	
	\begin{eqnarray*}
		M&=&\displaystyle\sup_{B(p,2R)}f^{-2},\nonumber\\
		P&=&\frac{(n+m)c_{1}^{2}}{4R^{2}\beta(1-\beta)}+\frac{(n-1+R\sqrt{nK})c_{1}+c_{2}+2c_{1}^{2}}{R^2},\nonumber\\
		Q&=&\frac{3\beta}{2} \left[\frac{n+m}{4}(k\gamma)^{4}\frac{(1-\beta)^{2}}{\beta^{4}}\varepsilon^{-1}\right]^{\frac{1}{3}}+\frac{(n+m)}{2}\beta(1-\varepsilon)^{-1}(1-\beta)^{-2}K^{2}, \nonumber\\
		S&=&\frac{\beta(n+m)}{2(1-\varepsilon)(1-\beta)^{2}}\Bigg{\{}\left[\frac{M\theta }{\beta}\left(\beta-1+\frac{2}{k}\right)\right]^{2}
		+\frac{2KM\theta}{\beta}\left(\beta-1+\frac{2}{k}\right)\Bigg{\}},\nonumber
	\end{eqnarray*}
$c_{1}$, $c_{2}$ are positive constants and $\varepsilon\in(0,1)$.
\end{theorem}

By letting $R\to\infty$ we get the following global gradient estimates.

\begin{cor}Let $(B^{n}\times_{f}F^{k},g, h,\rho)$ be a complete gradient almost Ricci soliton with nocompact base satisfying
	$$Ric_{g_B}+Hess_{g_B}h_{B}-\frac{1}{m}dh_{B}\otimes dh_{B}\geq -K,\qquad\Delta_{h}\rho_{B}\geq0,\qquad|\nabla\rho_{B}|\leq\gamma,$$
	for $K\geq0$ on $B$. Then, for any $\beta\in(0,1)$ with $\beta>1-\frac{2}{k}$, the warping function $f$ satisfies the following gradient estimates:
	
	\begin{itemize}
		\item If $\theta<0$, we have
		$$\beta\frac{|\nabla f|^{2}}{f^{2}}+\frac{\rho_B}{k}-\frac{\theta}{kf^2}\leq \sqrt{\frac{n+m}{2\beta k^{4}}}Q^{\frac{1}{2}},\qquad \text{in}\quad  B^{n},$$
		
		\item If $\theta=0$, we have
		$$\beta\frac{|\nabla f|^{2}}{f^{2}}+\frac{\rho_B}{k}\leq \sqrt{\frac{n+m}{2\beta k^{4}}}Q^{\frac{1}{2}},\qquad \text{in}\quad  B^{n},$$
		
		\item If $\theta>0$, assume that $f$ is bounded, then we have
			$$\beta\frac{|\nabla f|^{2}}{f^{2}}+\frac{\rho_B}{k}-\frac{\theta}{kf^2}\leq \frac{2(n+m)\theta M'}{k^{2}\beta}+\sqrt{\frac{n+m}{2\beta k^{4}}}(Q+S)^{\frac{1}{2}},\qquad \text{in}\quad  B^{n},$$
		
	\end{itemize}
 where
	
	\begin{eqnarray*}
		M'&=&\displaystyle\sup_{B}f^{-2},\nonumber\\
		Q&=&\frac{3\beta}{2} \left[\frac{n+m}{4}(k\gamma)^{4}\frac{(1-\beta)^{2}}{\beta^{4}}\varepsilon^{-1}\right]^{\frac{1}{3}}+\frac{(n+m)}{2}\beta(1-\varepsilon)^{-1}(1-\beta)^{-2}K^{2}, \nonumber\\
		S&=&\frac{\beta(n+m)}{2(1-\varepsilon)(1-\beta)^{2}}\Bigg{\{}\left[\frac{M\theta }{\beta}\left(\beta-1+\frac{2}{k}\right)\right]^{2}
		+\frac{2KM\theta}{\beta}\left(\beta-1+\frac{2}{k}\right)\Bigg{\}},\nonumber
	\end{eqnarray*}
and $\varepsilon\in(0,1)$.
\end{cor}

As an application we obtain the following two results.

\begin{cor}\label{cor}There is no complete gradient Ricci soliton $(B^{n}\times_{f}F^{k},g, h,\rho)$ satisfying $$ Ric_{g_B}+Hess_{g_B}h_B-\frac{1}{m}dh_B\otimes dh_B\geq 0, $$ if either $\rho=0$, $\theta<0$ or $\rho=\text{cte}>0$, $\theta=0$.
\end{cor}

\begin{example}\label{ex:1}
	Let $(\mathbb{H}^{n},g_{-1})$, $(\mathbb{R}^{n},g_{0})$ and $(\mathbb{S}^{n},g_{1})$ be, in this order, the hyperbolic  and Euclidean spaces as well as the round sphere, all with standard metrics. Consider the manifolds $(B,g_B)=(\mathbb{S}^{n}\times \mathbb{R}^n, g_1+g_{0})$, $(F,g_F)=(\mathbb{H}^k,g_{-1})$ and define
	\[h_v:\mathbb{S}^n\to\mathbb{R},\quad p\mapsto g_1(p,v),\]
	for any given $v\in\mathbb{S}^n$. Then, it's readily seen that
	\[h_{B}: \mathbb{S}^{n}\times \mathbb{R}^n\to\mathbb{R},\quad(x,y)\mapsto -\log{\left(\frac{3-h_{v}(x)}{n}\right)},\]
	satisfies 
	\[Ric_{g_B}+Hess_{g_B}h_B-dh_B\otimes dh_B\geq 0,\]
	and, therefore, by Corollary \ref{cor}, there cannot exist an $f:\mathbb{S}^n \times\mathbb{R}^n\to(0,\infty)$, not even a constant one, such that the warped product $ ((\mathbb{S}^n\times \mathbb{R}^n)\times_f\mathbb{H}^k,g) $ is a steady gradient Ricci soliton with potential $ h_B\circ\pi $.
\end{example}

Corollary \ref{cor} revels us that by starting with a quasi-Einstein manifold $(B,g_B)$ such that
\begin{equation*}
Ric_{g_B}+Hess_{g_B}h_B-\frac{1}{m}dh_B\otimes dh_B=0,
\end{equation*}
it is impossible to utilize such manifold as the base of a warped product gradient Ricci soliton $ (B^n\times F^k,g,h_B\circ\pi,\rho) $ with either a negative scalar curvature of the fiber manifold $ <0 $ and $ \rho=0 $ or null scalar curvature of the fiber manifold $ =0 $ and $ \rho>0 $.

In particular, for Einstein metrics, we foresee that


\begin{cor}
	\label{cor2}
	There is no Einstein warped product $(B^{n}\times_{f}F^{k},g)$ such that $\rho=0$ and $Ric_{g_B}\geq0$ with either $\theta<0$ or $ \theta=0 $ and $ f\neq constant $.
\end{cor}


\begin{example}
	By taking $ (B,g_B)=(\mathbb{R}^n,g_0) $ and $ (F,g_F)=(\mathbb{H}^k,g_{-1}) $ we see, now from Corollary \ref{cor2}, that none of the warped product manifolds $ (\mathbb{R}^n\times_f\mathbb{H}^k,g) $ are Ricci flat.
\end{example}

\section{Proofs}
\label{proofs}

\begin{myproof}{Theorem}{\ref{gradient}} Applying the change $f=v^{\frac{1}{k}}$ in equation \eqref{nazareno1} of Proposition \ref{lemma1} we obtain that
\begin{equation}\label{9}
\Delta_{h_{B}}v+\rho_{B} kv-\theta k v^{1-\frac{1}{k}}=0.
\end{equation}
Let $v$ a positive solution to \eqref{9}. Then, $u=\log v$ satisfies

\begin{equation*}
\Delta_{h_B}u=(\beta-1)|\nabla u|^{2}-L,\quad \text{where}\quad L:=\beta|\nabla u|^{2}+\rho_{B}k-\theta k e^{-\frac{2u}{k}},\quad \beta\in(0,1),\quad \beta>1-\frac{2}{k}.
\end{equation*}

Now, consider a cut-off function $\xi$ satisfying

\begin{equation*}
\xi(r)= \begin{cases} 1 &\text{if } r\in[0,1] \\[10pt]
0 & \text{if } r\in[2,\infty) 
\end{cases},\qquad -c_{1}\leq\frac{\xi'(r)}{\xi^{\frac{1}{2}}(r)}\leq0,\qquad -c_{2}\leq\xi''(r),\qquad c_{1},c_{2}\in (0,\infty)
\end{equation*}
and define
\begin{equation*}
\psi(x)=\xi\left(\frac{r(x)}{R}\right),
\end{equation*}
where $r(x)$ is the distance function from $p$ to $x$.
Using an argument of Calabi \cite{calabi1958extension} (see also Cheng and Yau \cite{cheng1975differential}), we can assume without loss of generality that the function $\psi$ is smooth in $B(p,2R)$. Then, the function defined by $G=\psi L$ is smooth in $B(p,2R)$. 

Let $x_{0}\in B(p,2R)$ be a point at which $G$ attains its maximum value $G_{max}$, and suppose that $G_{max}>0$ (otherwise the proof is trivial). At the point $x_{0}$, we have
\begin{equation*}
\nabla (G)=\psi\nabla L+L\nabla \psi=0.
\end{equation*}
Moreover,
\begin{equation}\label{general}
\begin{split}
0&\geq\Delta_{h_{B}}G,\\
&=\psi\Delta_{h_{B}}L+L\Delta_{h_{B}}\psi+2\langle\nabla \psi,\nabla L\rangle,\\
&=\psi\Delta_{h_B}L+L\Delta_{h_B}\psi-2L\frac{|\nabla\psi|^{2}}{\psi}.
\end{split}
\end{equation}
In order to estimate the right-hand side of \eqref{general} we prove the following lemma:

\begin{lema}\label{lemma2}Let $(B^{n},g_{B})$ be a complete noncompact Riemannian  manifold satisfying,
	
	\begin{equation}\label{m-Bakry}
	Ric_{g_B}+Hess_{g_B}h_B-\frac{1}{m}dh_{B}\otimes dh_{B}\geq -K,
	\end{equation}
	for $K\geq0$ in the metric ball $B(p,2R)\subset B$, and let $L$ and $\psi$ as above. Then, we have that

	\begin{equation}\label{a2'}
	\frac{|\nabla\psi|^{2}}{\psi}\leq\frac{c_{1}^{2}}{R^{2}},
	\end{equation}
	
	\begin{equation}\label{a3'}
	\Delta_{h_{B}}\psi\geq-\frac{(n-1+R\sqrt{nK})c_{1}+c_{2}}{R^2},
	\end{equation}
	
		\begin{equation}\label{a1}
		\begin{split}
		\Delta_{h_B}L\geq 2\beta\frac{(\Delta_{h_B}u)^{2}}{n+m}+2(1-\beta)k\langle\nabla u,\nabla \rho_{B}\rangle-2\langle\nabla u,\nabla L\rangle-2\beta K|\nabla u|^{2}
		+\\+k\Delta_{h_B}\rho_B-2\theta e^{-\frac{2u}{k}}\left(\big{(}\beta-1+\frac{2}{k}\big{)}|\nabla u|^{2}+L\right).
		\end{split}
	\end{equation}
\end{lema}

\begin{myproof}{Lemma}{\ref{lemma2}}
Equation \eqref{a2'} follows from the calculation

\begin{equation*}
\frac{|\nabla\psi|^{2}}{\psi}=\frac{1}{\xi}\left\langle \xi'\frac{\nabla r}{R}, \xi'\frac{\nabla r}{R}\right\rangle=\frac{(\xi')^2}{\xi}\frac{1}{R^2}\langle\nabla r,\nabla r\rangle\leq\frac{c_{1}^{2}}{R^{2}}.\\[10pt]
\end{equation*}

It has been shown by Qian \cite{qian1998comparison}, the following estimate

\begin{equation*}
\Delta_{h_{B}}r^{2}\leq n\left(1+\sqrt{1+\frac{4Kr^2}{n}}\right)
\end{equation*}
which implies
\begin{eqnarray}
\Delta_{h_B}r&=&\frac{1}{2r}\left(\Delta_{h_B}r^2-2|\nabla r|^{2}\right)\nonumber\\
&\leq&\frac{n-2}{2r}+\frac{n}{2r}\left(1+\sqrt{1+\frac{4Kr^2}{n}}\right)\nonumber\\
&=&\frac{n-1}{r}+\sqrt{nK}\nonumber.
\end{eqnarray}
Then, we obtain
\begin{equation*}
\Delta_{h_B}\psi=\frac{\xi''(r)|\nabla r|^{2}}{R^2}+\frac{\xi'(r)\Delta_{h_B}r}{R}\geq-\frac{(n-1+R\sqrt{nK})c_{1}+c_{2}}{R^2}.
\end{equation*}
which proves \eqref{a3'}.

From the Bochner formula for the $m$-Bakry-\'Emery Ricci tensor and
the lower bound hypothesis \eqref{m-Bakry}, we obtain

\begin{eqnarray*}
	\frac{1}{2}\Delta_{h_B}|\nabla u|^{2}\geq\frac{(\Delta_{h_B}u)^{2}}{n+m}+\langle\nabla u,\nabla\Delta_{h_B}u\rangle-K|\nabla u|^{2}.
\end{eqnarray*}
Therefore, 
\begin{equation*}
\begin{split}
\Delta_{h_B}L&=\beta\Delta_{h_B}|\nabla u|^{2}+k\Delta_{h_B}\rho_B-k\theta\Delta_{h_B}e^{-\frac{2u}{k}}\\[5pt]
&\geq 2\beta\frac{(\Delta_{h_B}u)^{2}}{n+m}+2\beta\langle\nabla u,\nabla\Delta_{h_B}u\rangle-2\beta K|\nabla u|^{2}+k\Delta_{h_B}\rho_B-k\theta\Delta_{h_B}e^{-\frac{2u}{k}}.
\end{split}
\end{equation*}
Notice that

\begin{eqnarray}
2\beta\langle\nabla u,\nabla\Delta_{h_B}u\rangle&=&2\beta\langle\nabla u,\nabla\left[\left(1-\frac{1}{\beta}\right)\left(-k\rho_B+k\theta e^{-\frac{2u}{k}}\right)-\frac{L}{\beta}\right]\rangle\nonumber\\[5pt]
&=&2(1-\beta)k\langle\nabla u,\nabla \rho_{B}\rangle+2(\beta-1) k\theta\langle\nabla u,\nabla e^{-\frac{2u}{k}}\rangle-2\langle\nabla u,\nabla L\rangle\nonumber\\[5pt]
&=&2(1-\beta)k\langle\nabla u,\nabla \rho_{B}\rangle-4(\beta-1) \theta e^{-\frac{2u}{k}}|\nabla u|^{2}-2\langle\nabla u,\nabla F\rangle,\nonumber
\end{eqnarray}
and

\begin{eqnarray}
\Delta_{h_B}e^{-\frac{2u}{k}}&=&\frac{2}{k}e^{-\frac{2u}{k}}\left[\frac{2}{k}|\nabla u|^{2}-\Delta_{h_B}u\right]\nonumber\\
&=&\frac{2}{k}e^{-\frac{2u}{k}}\left[\frac{2}{k}|\nabla u|^{2}-(\beta-1)|\nabla u|^{2}+L\right]\nonumber\\
&=&\frac{2}{k}e^{-\frac{2u}{k}}\left[\left(\frac{2}{k}-\beta+1\right)|\nabla u|^{2}+L\right].\nonumber
\end{eqnarray}
It follows that

\begin{equation*}
\begin{split}
\Delta_{h_B}L\geq 2\beta\frac{(\Delta_{h_B}u)^{2}}{n+m}+2(1-\beta)k\langle\nabla u,\nabla \rho_{B}\rangle-2\langle\nabla u,\nabla L\rangle-2\beta K|\nabla u|^{2}
+\\+k\Delta_{h_B}\rho_B-2\theta e^{-\frac{2u}{k}}\left(\big{(}\beta-1+\frac{2}{k}\big{)}|\nabla u|^{2}+L\right),
\end{split}
\end{equation*}
which completes the proof of lemma.
\end{myproof}\\[1pt]

Continuing, by Lemma \ref{lemma2} and \eqref{general}, we obtain at the point $x_{0}$,

\begin{eqnarray*}
\psi\Bigg{\{}2\beta\frac{(\Delta_{h_B}u)^{2}}{n+m}+2(1-\beta)k\langle\nabla u,\nabla \rho_{B}\rangle-2\langle\nabla u,\nabla L\rangle-2\beta K|\nabla u|^{2}+k\Delta_{h_B}\rho_B+\\-2\theta e^{-\frac{2u}{k}}\left[\big{(}\beta-1+\frac{2}{k}\big{)}|\nabla u|^{2}+L\right]\Bigg{\}}\leq LH,
\end{eqnarray*}
where

\begin{equation*}
H=\left(\frac{(n-1+R\sqrt{nK})c_{1}+c_{2}+2c_{1}^{2}}{R^2}\right).
\end{equation*}

From the fact that $0\leq\psi\leq1$, we have
$$
-2\psi\langle\nabla u,\nabla L\rangle=2L\langle\nabla u, \nabla \psi\rangle\geq-2L|\nabla u||\nabla\psi|\geq-\frac{2c_{1}}{R}\psi^{\frac{1}{2}}L|\nabla u|.
$$
Then

\begin{equation}\label{final}
\begin{split}
	2\beta\psi\frac{(\Delta_{h_B}u)^{2}}{n+m}+2(1-\beta)k\psi\langle\nabla u,\nabla \rho_{B}\rangle
	-\frac{2c_{1}}{R}\psi^{\frac{1}{2}}L|\nabla u|
	-2\beta\psi K|\nabla u|^{2}+\\
	+k\psi\Delta_{h_B}\rho_B-2\theta e^{-\frac{2u}{k}}\psi\left[\big{(}\beta-1+\frac{2}{k}\big{)}|\nabla u|^{2}+L\right]\leq LH.
\end{split}
\end{equation}
In the sequel we distinguish between two cases: (a) $\theta<0$ and (b) $\theta\geq0$.

Case (a): $\theta<0$. Since

$$\Delta_{h_B}\rho_{B}\geq0,\qquad|\nabla \rho_{B}|\leq\gamma(2R),$$
then \eqref{final} yields

\begin{eqnarray*}
		2\beta\psi\frac{(\Delta_{h_B}u)^{2}}{n+m}-2(1-\beta)k\psi\gamma|\nabla u|
		-\frac{2c_{1}}{R}\psi^{\frac{1}{2}}L|\nabla u|
		-2\beta\psi K|\nabla u|^{2}\leq LH.
\end{eqnarray*}
Multiplying both sides of the above equation by $\psi$ and using the fact that $0\leq\psi\leq 1$, we obtain

\begin{eqnarray*}
	2\beta\frac{(\psi\Delta_{h_B}u)^{2}}{n+m}+2(\beta-1)k\psi^{\frac{1}{2}}\gamma|\nabla u|
	-\frac{2c_{1}}{R}\psi^{\frac{3}{2}}L|\nabla u|
	-2\beta\psi K|\nabla u|^{2}\leq \psi LH.
\end{eqnarray*}
Let

$$y=\psi|\nabla u|^{2},\qquad z=\psi(-k\rho_B+k\theta e^{-\frac{2u}{k}}).$$
Then we have

\begin{eqnarray*}
	\frac{2\beta}{n+m}\Big\{ (y-z)^{2}+\frac{(\beta-1)k\gamma(n+m)y^{\frac{1}{2}}}{\beta}-\frac{(n+m)c_{1}}{R}y^{\frac{1}{2}}\left(y-\frac{z}{\beta}\right)-(n+m)Ky\Big\}\leq \psi LH.
\end{eqnarray*}
From Li-Yau's arguments (\cite{li1986parabolic}, pg.161-162), we know that

\begin{eqnarray*}
(y-z)^{2}-(n+m)c_{1}R^{-1}y^{\frac{1}{2}}\left(y-\frac{z}{\beta}\right)-(n+m)Ky-(n+m)\left(\frac{1}{\beta}-1\right)\gamma y^{\frac{1}{2}}\geq\nonumber\\
\geq\left(\frac{1}{\beta}\right)^{-2}\left(y-\frac{z}{\beta}\right)^{2}-\frac{(n+m)^{2}}{8}c_{1}^{2}\left(\frac{1}{\beta}\right)^{2}\left(\frac{1}{\beta}-1\right)^{-1}R^{-2}\left(y-\frac{z}{\beta}\right)+\nonumber\\
-\frac{3}{4}4^{-\frac{1}{3}}(n+m)^{\frac{4}{3}}\left[\gamma^{4}\left(\frac{1}{\beta}-1\right)^{2}\left(\frac{1}{\beta}\right)^{2}\varepsilon^{-1}\right]^{\frac{1}{3}}-\frac{(n+m)^{2}}{4}(1-\varepsilon)^{-1}\left(\frac{1}{\beta}-1\right)^{-2}\left(\frac{1}{\beta}\right)^{2}K^{2}.
\end{eqnarray*}
for any $0<\varepsilon<1$.

Therefore in our context, we obtain
\begin{eqnarray*}
	\frac{2\beta}{n+m}\Bigg\{(\beta y-z)^{2}-\frac{(n+m)^{2}c_{1}^{2}}{8R^{2}\beta^{2}(1-\beta)}(\beta y-z)-\frac{3}{4}4^{-\frac{1}{3}}(n+m)^{\frac{4}{3}}\left[\left(k\gamma\right)^{4}\frac{(1-\beta)^{2}}{\beta^{4}}\varepsilon^{-1}\right]^{\frac{1}{3}}+\nonumber\\
	-\frac{(n+m)^{2}}{4}(1-\varepsilon)^{-1}(1-\beta)^{-2}K^{2}\Bigg\}
	\leq (\beta y-z)H.
\end{eqnarray*}
Hence,

\begin{equation*}
\frac{2\beta}{n+m}(\psi L)^{2}-P(\psi L)-Q\leq0,
\end{equation*}
where,

\begin{eqnarray*}
P&=&\frac{(n+m)c_{1}^{2}}{4R^{2}\beta(1-\beta)}+H,\nonumber\\
Q&=&\frac{3\beta}{2} \left[\frac{n+m}{4}\left(k\gamma\right)^{4}\frac{(1-\beta)^{2}}{\beta^{4}}\varepsilon^{-1}\right]^{\frac{1}{3}}+\frac{\beta(n+m)}{2}(1-\varepsilon)^{-1}(1-\beta)^{-2}K^{2}. \nonumber
\end{eqnarray*}

Using the inequality $az^{2}-bz\leq c$, one obtain $z\leq\frac{2b}{a}+\sqrt{\frac{c}{a}}$. Then

\begin{eqnarray*}
\sup_{x\in B(p,R)}L(x)\leq(\psi L)(x_{0})\leq \frac{n+m}{\beta}P+\sqrt{\frac{n+m}{2\beta}}Q^{\frac{1}{2}},
\end{eqnarray*}
and hence,

$$\beta|\nabla u|^{2}+k\rho_B-k\theta e^{-\frac{2u}{k}}\leq \frac{n+m}{\beta}P+\sqrt{\frac{n+m}{2\beta}}Q^{\frac{1}{2}}.$$

Replacing the function $u=\log f^{\frac{d+1}{2}}$ back into the above equation we obtain the desired estimate when $\theta<0$.

Case (b): $\theta\geq0$. Since
$$\Delta_{h_B}\rho_{B}\geq0,\qquad|\nabla \rho_{B}|\leq\gamma(2R),$$
then \eqref{final} yields

\begin{equation*}
\begin{split}
2\beta\psi\frac{(\Delta_{h_B}u)^{2}}{n+m}-2(1-\beta)k\psi|\nabla u|\gamma
-\frac{2c_{1}}{R}\psi^{\frac{1}{2}}L|\nabla u|
-2\beta\psi K|\nabla u|^{2}+\\
-2\theta e^{-\frac{2u}{k}}\psi\left[\big{(}\beta-1+\frac{2}{k}\big{)}|\nabla u|^{2}+L\right]\leq LH.
\end{split}
\end{equation*}
Multiplying both sides of the above equation by $\psi$, and using the fact that $0\leq\psi\leq 1$, we obtain

\begin{eqnarray*}
	2\beta\frac{(\psi\Delta_{h_B}u)^{2}}{n+m}+2(\beta-1)k\psi^{\frac{1}{2}}\gamma|\nabla u|-\frac{2c_{1}}{R}\psi^{\frac{3}{2}}L|\nabla u|-2\beta\psi K|\nabla u|^{2}+\\
	-2\theta M\psi\left[\big{(}\beta-1+\frac{2}{k}\big{)}|\nabla u|^{2}+L\right]\leq \psi LH,
\end{eqnarray*}
where $M=\sup_{B(p,2R)} e^{-\frac{2u}{k}}$.

Let
$$y=\psi|\nabla u|^{2},\qquad z=\psi(-k\rho_B+k\theta e^{-\frac{2u}{k}}).$$
Then we have

\begin{eqnarray*}
	\frac{2\beta}{n+m}\Bigg\{(y-z)^{2}-\frac{c_{1}(n+m)}{R}y^{\frac{1}{2}}\left(y-\frac{z}{\beta}\right)
	+(n+m)\frac{(\beta-1)}{\beta}k\gamma y^{\frac{1}{2}}+\nonumber\\
	-(n+m)\left[\frac{\theta M}{\beta}(\beta-1+\frac{2}{k})+K\right]y
	\Bigg\}\leq \psi L(H+2\theta M).
\end{eqnarray*}
Hence, again by the Li-Yau arguments, we get that

\begin{eqnarray*}
	\frac{2\beta}{n+m}\Bigg\{(\beta y-z)^{2}-\frac{(n+m)^{2}c_{1}^{2}}{8R^{2}\beta^{2}(1-\beta)}(\beta y-z)-\frac{3}{4}4^{-\frac{1}{3}}(n+m)^{\frac{4}{3}}\left[\left(k\gamma\right)^{4}\frac{(1-\beta)^{2}}{\beta^{4}}\varepsilon^{-1}\right]^{\frac{1}{3}}+\nonumber\\
	-\frac{(n+m)^{2}}{4}(1-\varepsilon)^{-1}(1-\beta)^{-2}\left[\frac{\theta M}{\beta}(\beta-1+\frac{2}{k})+K\right]^{2}\Bigg\}
	\leq (\beta y-z)H,
\end{eqnarray*}
which means that 

\begin{eqnarray*}
	\frac{2\beta}{n+m}(\psi L)^{2}-\big\{P+2\theta M\big\}(\psi L)-(Q+S)\leq 0,
\end{eqnarray*}
where

\[S=\frac{\beta(n+m)}{2(1-\varepsilon)(1-\beta)^{2}}\Bigg{\{}\left[\frac{M\theta }{\beta}\left(\beta-1+\frac{2}{k}\right)\right]^{2}
+\frac{2KM\theta}{\beta}\left(\beta-1+\frac{2}{k}\right)\Bigg{\}}.\]
Therefore

\begin{eqnarray*}
	\beta|\nabla u|^{2}+k\rho_B-k\theta e^{-\frac{2u}{k}}\leq \frac{n+m}{\beta}(P+2\theta M)+\sqrt{\frac{n+m}{2\beta}}(Q+S)^{\frac{1}{2}}.
\end{eqnarray*}
Replacing the function $u=\log f^{\frac{d+1}{2}}$ back into the above equation we obtain the desired estimate when $\theta\geq0$.
\end{myproof}

\end{document}